\documentclass[12pt,twoside]{article}
\usepackage[english]{babel}
\usepackage[latin1]{inputenc}
\usepackage{amsmath}
\usepackage{graphicx}
\usepackage{times,amssymb,amscd}
\newtheorem{thm}{Theorem}[section]

\newtheorem{lem}[thm]{Lemma}

\numberwithin{equation}{section}

\DeclareMathOperator{\arccosh}{arccosh}

\newcommand{\bA}{\mathbf{A}}

\newcommand{\bE}{\mathbf{E}}

\newcommand{\bH}{\mathbf{H}}

\newcommand{\bL}{\mathbf{L}}

\newcommand{\bR}{\mathbf{R}}
\newcommand{\bS}{\mathbf{S}}
\newcommand{\bV}{\mathbf{V}}

\newcommand{\be}{\mathbf{e}}

\newcommand{\bx}{\mathbf{x}}
\newcommand{\by}{\mathbf{y}}

\newcommand{\bT}{\mathbf{T}}

\newcommand{\BV}{\boldsymbol{V}}

\newcommand{\Be}{\boldsymbol{e}}
\newcommand{\Bu}{\boldsymbol{u}}
\newcommand{\Bv}{\boldsymbol{v}}

\newcommand{\cA}{\mathcal{A}}

\newcommand{\cP}{\mathcal{P}}
\newcommand{\cS}{\mathcal{S}}

\newcommand{\EUC}{\mathbf E^3}
\newcommand{\SPH}{\bS^3}
\newcommand{\HYP}{\bH^3}
\newcommand{\SXR}{\bS^2\!\times\!\bR}
\newcommand{\HXR}{\bH^2\!\times\!\bR}
\newcommand{\SLR}{\widetilde{\bS\bL_2\bR}}
\newcommand{\NIL}{\mathbf{Nil}}
\newcommand{\SOL}{\mathbf{Sol}}

\begin{document}
\pagestyle{myheadings}
\markboth{\centerline{Jen\H o Szirmai}}
{Apollonius surfaces, circumscribed spheres of tetrahedra, Menelaus' and Ceva's theorems $\dots$}
\title
{Apollonius surfaces, circumscribed spheres of tetrahedra, Menelaus' and Ceva's theorems in $\SXR$ and $\HXR$ geometries
}

\author{Jen\H o Szirmai \\
\normalsize Budapest University of Technology and \\
\normalsize Economics Institute of Mathematics, \\
\normalsize Department of Geometry \\
\normalsize Budapest, P. O. Box: 91, H-1521 \\
\normalsize szirmai@math.bme.hu
\date{\normalsize{\today}}}
\maketitle
\begin{abstract}
In the present paper we study $\SXR$ and $\HXR$ geometries, which are homogeneous Thurston 3-geometries.
We define and determine the generalized Apollonius surfaces and with them define the ``surface of a geodesic triangle". 
Using the above Apollonius surfaces we develop a procedure to determine the centre and the radius of the circumscribed geodesic sphere of an arbitrary $\SXR$ and $\HXR$ tetrahedron.
Moreover, we generalize the famous Menelaus' and Ceva's theorems for geodesic triangles in both spaces. 
In our work we will use the projective model of $\SXR$ and $\HXR$ geometries described by E. Moln\'ar in \cite{M97}.

\end{abstract}
\newtheorem{Theorem}{Theorem}[section]
\newtheorem{corollary}[Theorem]{Corollary}
\newtheorem{lemma}[Theorem]{Lemma}
\newtheorem{exmple}[Theorem]{Example}
\newtheorem{definition}[Theorem]{Definition}
\newtheorem{rmrk}[Theorem]{Remark}
\newtheorem{proposition}[Theorem]{Proposition}
\newenvironment{remark}{\begin{rmrk}\normalfont}{\end{rmrk}}
\newenvironment{example}{\begin{exmple}\normalfont}{\end{exmple}}
\newenvironment{acknowledgement}{Acknowledgement}

\section{Introduction}
\label{section1}
The classical definition of the Apollonius circle in the Euclidean plane $\mathbf{E}^2$ is  
the set of all points of $\mathbf{E}^2$ whose distances from two fixed points are in a constant ratio $\lambda\in\mathbf{R}^+$. This definition can be extended
in a natural way to the Thurston geometries 
$$
\EUC,\SPH,\HYP,\SXR,\HXR,\NIL,\SLR,\SOL.
$$
\begin{definition}
The Apollonius surface in the Thurston geometry $X$ is  
the set of all points of $X$ whose geodesic distances from two fixed points are in a constant ratio $\lambda\in\mathbf{R}^+$.
\end{definition}
\begin{rmrk}
A special case of Apollonius surfaces is the geodesic-like bisector (or equidistant) surface ($\lambda=1)$
of two arbitrary points of $X$. These surfaces have an important role in structure of Dirichlet - Voronoi (briefly, D-V) cells.

The D-V-cells
are relevant in the study of tilings, ball packing and ball covering. E.g. if the point set is the orbit of a point - generated by
a discrete isometry group of $X$ - then we obtain a monohedral D-V cell decomposition (tiling) of the considered space and it is interesting to examine its
optimal ball packing and covering. In $3$-dimensional spaces of constant curvature, the
D-V cells have been widely investigated, but in the other Thurston geometries $\SXR$, $\HXR$,
$\NIL$, $\SOL$, $\SLR$ there are few results on this topic. 

In \cite{PSSz10}, \cite{PSSz11-1}, \cite{PSSz11-2} we studied the geodesic-like 
equidistant surfaces in $\SXR$, $\HXR$ and $\NIL$ geometries, and in \cite{Sz19}, \cite{VSz19}
the translation-like equidistant surfaces in $\SOL$ and $\NIL$ geometries. 
\end{rmrk}

In the present paper, we are interested in {Apollonius surfaces, \it geodesic triangles and their surfaces, generalized Menelaus' and Ceva's theorems} in $\SXR$ and $\HXR$ spaces \cite{S,T}. 

In Section 2 we describe the projective model and the isometry group of the considered geometries,
moreover, we give an overview about its geodesic curves.
In Section 3 we study the generallized Apollonius surfaces and their properties in the considered spaces and using them we define the surfaces of geodesic triangles.

In Section 4 we generalize and prove the theorems of Menelaus and Ceva for geodesic triangles in $\SXR$ and $\HXR$ spaces.

The computation and the proof is based on the projective model 
of $\SXR$ and $\HXR$ geometries described by E. Moln\'ar in \cite{M97}.

\section{Projective model of $\HXR$ and $\SXR$ spaces}
E. {Moln\'ar} has shown in \cite{M97}, that the homogeneous 3-spaces
have a unified interpretation in the projective 3-sphere $\mathcal{PS}^3(\bV^4,\BV_4, \mathbf{R})$. 
In our work we shall use this projective model of $\SXR$ and $\HXR$ geometries. 
The Cartesian homogeneous coordinate simplex is given by $E_0(\be_0)$,$E_1^{\infty}(\be_1)$,$E_2^{\infty}(\be_2)$,
$E_3^{\infty}(\be_3)$, $(\{\be_i\}\subset \bV^4$ \ and $\text{with the unit point}$ $E(\be = \be_0 + \be_1 + \be_2 + \be_3 ))$. 
Moreover, $\by=c\bx$ with $0<c\in \mathbf{R}$ (or $c\in\mathbf{R}\setminus\{0\})$
defines a point $(\bx)=(\by)$ of the projective 3-sphere $\cP \cS^3$ (or that of the projective space $\cP^3$ where opposite rays
$(\bx)$ and $(-\bx)$ are identified). 
The dual system $\{(\Be^i)\}\subset \BV_4$ describes the simplex planes, especially the plane at infinity 
$(\Be^0)=E_1^{\infty}E_2^{\infty}E_3^{\infty}$, and generally, $\Bv=\Bu\frac{1}{c}$ defines a plane $(\Bu)=(\Bv)$ of $\cP \cS^3$
(or that of $\cP^3$). Thus $0=\bx\Bu=\by\Bv$ defines the incidence of point $(\bx)=(\by)$ and plane
$(\Bu)=(\Bv)$, as $(\bx) \text{I} (\Bu)$ also denotes it. Thus {$\SXR$} can be visualized in the affine 3-space $\bA^3$
(so in $\bE^3$) as well.
\subsection{Geodesic curves and spheres in $\SXR$ space}
In this section we recall the important notions and results from the papers \cite{M97}, \cite{PSSz10}, \cite{Sz13-1}, \cite{Sz11-1}, \cite{Sz11-2}. 

The well-known infinitezimal arc-length square at any point of $\SXR$ as follows
\begin{equation}
   \begin{gathered}
     (ds)^2=\frac{(dx)^2+(dy)^2+(dz)^2}{x^2+y^2+z^2}.
       \end{gathered} \tag{2.1}
     \end{equation}
We shall apply the usual geographical coordiantes $(\phi, \theta), ~ (-\pi < \phi \le \pi, ~ -\frac{\pi}{2}\le \theta \le \frac{\pi}{2})$ 
of the sphere with the fibre coordinate $t \in \bR$. We describe points in the above coordinate system in our model by the following equations: 
\begin{equation}
x^0=1, \ \ x^1=e^t \cos{\phi} \cos{\theta},  \ \ x^2=e^t \sin{\phi} \cos{\theta},  \ \ x^3=e^t \sin{\theta} \tag{2.2}.
\end{equation}
Then we have $x=\frac{x^1}{x^0}=x^1$, $y=\frac{x^2}{x^0}=x^2$, $z=\frac{x^3}{x^0}=x^3$, i.e. the usual Cartesian coordinates.
We obtain by \cite{M97} that in this parametrization the infinitezimal arc-length square 
at any point of $\SXR$ is the following
\begin{equation}
   \begin{gathered}
      (ds)^2=(dt)^2+(d\phi)^2 \cos^2 \theta +(d\theta)^2.
       \end{gathered} \tag{2.3}
     \end{equation}
The geodesic curves of $\SXR$ are generally defined as having locally minimal arc length between their any two (near enough) points. 
The equation systems of the parametrized geodesic curves $\gamma(t(\tau),\phi(\tau),\theta(\tau))$ in our model can be determined by the 
general theory of Riemann geometry (see \cite{KN}, \cite{Sz11-2}).

Then by (2.1-2) we get the equation systems of a geodesic curve in our Euclidean model (see \cite{Sz11-1}):
\begin{equation}
  \begin{gathered}
   x(\tau)=e^{\tau \sin{v}} \cos{(\tau \cos{v})}, \\ 
   y(\tau)=e^{\tau \sin{v}} \sin{(\tau \cos{v})} \cos{u}, \\
   z(\tau)=e^{\tau \sin{v}} \sin{(\tau \cos{v})} \sin{u},\\
   -\pi < u \le \pi,\ \ -\frac{\pi}{2}\le v \le \frac{\pi}{2}. \tag{2.4}
  \end{gathered}
\end{equation}
\begin{definition}
The distance $d^{\SXR}(P_1,P_2)$ between the points $P_1$ and $P_2$ is defined by the arc length of the geodesic curve 
from $P_1$ to $P_2$.
\end{definition}
\begin{definition}
 The geodesic sphere of radius $\rho$ (denoted by $S_{P_1}(\rho)$) with center at the point $P_1$ is defined as the set of all points 
 $P_2$ in the space with the condition $d^{\SXR}(P_1,P_2)=\rho$. Moreover, we require that the geodesic sphere is a simply connected 
 surface without selfintersection in $\SXR$.
 \end{definition}
 \begin{definition}
 The body of the geodesic sphere of centre $P_1$ and of radius $\rho$ in the $\SXR$ space is called geodesic ball, denoted by $B_{P_1}(\rho)$,
 i.e., $Q \in B_{P_1}(\rho)$ iff $0 \leq d(P_1,Q) \leq \rho$.
 \end{definition}
\begin{proposition}
The geodesic sphere and ball of radius $\rho$ exists in the $\SXR$ space if and only if $\rho \in [0,\pi].$
\end{proposition}
\subsection{Geodesic curves and spheres of $\HXR$ geometry}
In this section we recall the important notions and results from the papers \cite{M97},  \cite{PSSz11-2}, \cite{Sz12-1}.

The points of $\HXR$ space, forming an open cone solid in the projective space $\mathcal{P}^3$, are the following:
\begin{equation}
\HXR:=\big\{ X(\bx=x^i \be_i)\in \mathcal{P}^3: -(x^1)^2+(x^2)^2+(x^3)^2<0<x^0,~x^1 \big\}. \notag
\end{equation}
In this context E. Moln\'ar \cite{M97} has derived the infinitezimal arc-length square at any point of $\HXR$ as follows
\begin{equation}
   \begin{gathered}
     (ds)^2=\frac{1}{(-x^2+y^2+z^2)^2}\cdot [(x)^2+(y)^2+(z)^2](dx)^2+ \\ + 2dxdy(-2xy)+2dxdz (-2xz)+ [(x)^2+(y)^2-(z)^2] (dy)^2+ \\ 
     +2dydz(2yz)+ [(x)^2-(y)^2+(z)^2](dz)^2.
       \end{gathered} \tag{2.5}
     \end{equation}
This becomes simpler in the following special (cylindrical) coordinates $(t, r, \alpha)$, $(r \ge 0, ~ -\pi < \alpha \le \pi)$ 
with the fibre coordinate $t \in \bR$. We describe points in our model by the following equations: 
\begin{equation}
x^0=1, \ \ x^1=e^t \cosh{r},  \ \ x^2=e^t \sinh{r} \cos{\alpha},  \ \ x^3=e^t \sinh{r} \sin{\alpha}  \tag{2.6}.
\end{equation}
Then we have $x=\frac{x^1}{x^0}=x^1$, $y=\frac{x^2}{x^0}=x^2$, $z=\frac{x^3}{x^0}=x^3$, i.e. the usual Cartesian coordinates.
We obtain by \cite{M97} that in this parametrization the infinitezimal arc-length square by (2.5)
at any point of $\HXR$ is the following
\begin{equation}
   \begin{gathered}
      (ds)^2=(dt)^2+(dr)^2 +\sinh^2{r}(d\alpha)^2.
       \end{gathered} \tag{2.7}
     \end{equation}
The geodesic curves of $\HXR$ are generally defined as having locally minimal arc length between their any two (near enough) points. 
The equation systems of the parametrized geodesic curves $\gamma(t(\tau),r(\tau),\alpha(\tau))$ in our model can be determined by the 
general theory of Riemann geometry (see \cite{Sz12-1}):

Then by (2.6-7) we get the equation systems of a geodesic curve in our model \cite{Sz12-1}:
\begin{equation}
  \begin{gathered}
   x(\tau)=e^{\tau \sin{v}} \cosh{(\tau \cos{v})}, \\ 
   y(\tau)=e^{\tau \sin{v}} \sinh{(\tau \cos{v})} \cos{u}, \\
   z(\tau)=e^{\tau \sin{v}} \sinh{(\tau \cos{v})} \sin{u},\\
   -\pi < u \le \pi,\ \ -\frac{\pi}{2}\le v \le \frac{\pi}{2}. \tag{2.8}
  \end{gathered}
\end{equation}
\begin{definition}
The distance $d^{\HXR}(P_1,P_2)$ between the points $P_1$ and $P_2$ is defined by the arc length of the geodesic curve 
from $P_1$ to $P_2$.
\end{definition}
\begin{definition}
 The geodesic sphere of radius $\rho$ (denoted by $S^{\HXR}_{P_1}(\rho)$) with centre at the point $P_1$ is defined as the set of all points 
 $P_2$ in the space with the condition $d^{\HXR}(P_1,P_2)=\rho$. Moreover, we require that the geodesic sphere is a simply connected 
 surface without selfintersection in $\HXR$ space.
 \end{definition}
 \begin{rmrk}
 In this paper we consider only the usual spheres with "proper centre", i.e. $P_1 \in \HXR$. 
 If the centre of a "sphere" lie on the absolute quadric or lie out of our model the notion of the "sphere" (similarly to the hyperbolic space),
 can be defined, but that cases we shall study in a forthcoming work.
 \end{rmrk}
 \begin{definition}
 The body of the geodesic sphere of centre $P_1$ and of radius $\rho$ in $\HXR$ space is called geodesic ball, denoted by $B_{P_1}(\rho)$,
 i.e. $Q \in B_{P_1}(\rho)$ iff $0 \leq d(P_1,Q) \leq \rho$.
 \end{definition}
 \begin{proposition}
 $S(\rho)$ is a simply connected surface in $\mathbf{E}^3$ for $\rho >0$.
 \end{proposition}
\begin{rmrk}
$\SXR$ and $\HXR$ are affine metric spaces (affine-projective spaces -- in the sense of the unified formulation of \cite{M97}). Therefore their linear, affine, unimodular,
etc. transformations are defined as those of the embedding affine space.
\end{rmrk}
\section{Apollonius surfaces in $\SXR$ and $\HXR$ geometries}
The generalization of the classical definition of the Apollonius circle of the Euclidean plane $\mathbf{E}^2$ to Thurston geometries is the following  
\begin{definition}
The Apollonius surface $\mathcal{A}\cS^X_{P_1P_2}(\lambda)$ in the Thurston geometry $X$ is  
the set of all points of $X$ whose geodesic distances from two fixed points are in a constant ratio $\lambda\in\mathbf{R}^+_0$ where $X\in \EUC,\SPH,\HYP,\SXR,\HXR,\NIL,\SLR,\SOL.
$ i.e. $\mathcal{A}\cS^X_{P_1P_2}(\lambda)$ of two arbitrary points $P_1,P_2 \in X$ consists of all points $P'\in X$,
for which $d^X(P_1,P')=\lambda \cdot d^X(P',P_2)$ ($\lambda\in [0,\infty$) where $d^X$ is the corresponding distance function of $X$. If $\lambda=0$, 
then $\mathcal{A}\cS^X_{P_1P_2}(0):=P_1$ and it is clear, that in case $\lambda \to \infty$ then $d(P',P_2) \to 0$ therefore we say $\mathcal{A}\cS^X_{P_1P_2}(\infty):=P_2$.
\end{definition}

Now, we consider only the $\SXR$ and $\HXR$ geometries and can be assumed by the homogeneity of the above geometries that the 
starting point of a given geodesic curve segment in both geometries is $P_1=(1,1,0,0)$.
The other endpoint will be given by its homogeneous coordinates $P_2=(1,x,y,z)$. In order to obtain the Apollonius surfaces of two given points $P_1,P_2$ to a given costant ratio
$\lambda\in\mathbf{R}^+$ we consider a geodesic curve segment $g^{X}_{P_1P_2}$ ($ X \in \{ \SXR, \HXR \})$ and determine its
parameters $(u,v,\tau)$ expressed by the real coordinates 
$x$, $y$, $z$ of $P_2$. 

We obtain directly by equation system (2.4) and (2.8) the following lemmas (see PSSz11-2,PSSz10,\cite{Sz12-1,Sz20}):
\begin{lem}
Let $(1,x,y,z)$ $(x,y,z \in \bR, x^2+y^2+z^2 \ne 0)$ be the homogeneous coordinates of the point $P \in \SXR$. The paramerters of the 
corresponding geodesic curve $g^{\SXR}{(A_1,P)}$ are the following:
\begin{enumerate}
\item $y,z \in \bR \setminus \{0\}$ and $x^2+y^2+z^2 \ne 1$;
\begin{equation}
\begin{gathered}
v=\mathrm{arctan}\Big(\frac{\log \sqrt{x^2+y^2+z^2}}{\mathrm{arccos}\frac{x}{\sqrt{x^2+y^2+z^2}}}\Big),~u=\mathrm{arctan}\Big(\frac{z}{y}\Big),\\
\tau=\frac{\log \sqrt{x^2+y^2+z^2}}{\sin v}, ~ \text{where} ~ -\pi < u \le \pi, ~ -\pi/2\le v \le \pi/2, ~ \tau \in \bR^+.
\end{gathered} \tag{3.1}
\end{equation} 
\item $y=0$, $z\ne 0$ and $x^2+z^2 \ne 1$;
\begin{equation}
\begin{gathered}
u=\frac{\pi}{2}, ~v=\mathrm{arctan}\Big(\frac{\log \sqrt{x^2+z^2}}{\mathrm{arccos}\frac{x}{\sqrt{x^2+z^2}}}\Big),\\
\tau=\frac{\log \sqrt{x^2+z^2}}{\sin v}, ~ \text{where}  ~ -\pi/2\le v \le \pi/2, ~ \tau \in \bR^+.
\end{gathered} \tag{3.2}
\end{equation}
\item $y=0$, $z\ne 0$ and $x^2+z^2 = 1$;
\begin{equation}
\begin{gathered}
u=\frac{\pi}{2}, ~v=0, ~ \tau=\arccos(x), ~ \tau \in \bR^+.
\end{gathered} \tag{3.3}
\end{equation}
\item $y, z =0$;
\begin{equation}
u=0,~v=\frac{\pi}{2},~\tau=\log \sqrt{x^2+y^2+z^2}, ~ \tau \in \bR^+. \tag{3.4} 
\end{equation}
\item $x=0,~y=0$ and $z \ne 1$;
\begin{equation}
\begin{gathered}
u=\frac{\pi}{2},~ v=\mathrm{arctan}\frac{2  \log|z|}{\pi}, ~\tau=\frac{\log|z|}{\sin v}, \\ -\pi/2\le v \le \pi/2, ~ \tau \in \bR^+. \tag{3.5}
\end{gathered}
\end{equation}
\end{enumerate} ~ ~ $\square$
\end{lem}
\begin{lem}
Let $(1,x,y,z)$ $(x,y,z \in \bR, x^2-y^2-z^2 \ge 0,~ x \ge 0)$ be the homogeneous coordinates of the point $P \in \HXR$. The paramerters of the 
corresponding geodesic curve $g^{\HXR}{(A_1,P)}$ are the following:
\begin{enumerate}
\item $y,z \in \bR \setminus \{0\}$ and $x^2-y^2-z^2 \ne 1$;
\begin{equation}
\begin{gathered}
v=\mathrm{arctan}\Big(\frac{\log \sqrt{x^2-y^2-z^2}}{\mathrm{arccosh}\frac{x}{\sqrt{x^2-y^2-z^2}}}\Big),~u=\mathrm{arctan}\Big(\frac{z}{y}\Big),\\
\tau=\frac{\log \sqrt{x^2-y^2-z^2}}{\sin v}, ~ \text{where} ~ -\pi < u \le \pi, ~ -\pi/2\le v \le \pi/2, ~ \tau \in \bR^+.
\end{gathered} \tag{3.6}
\end{equation} 
\item $y=0$, $z\ne 0$ and $x^2-z^2 \ne 1$;
\begin{equation}
\begin{gathered}
u=\frac{\pi}{2}, ~v=\mathrm{arctan}\Big(\frac{\log \sqrt{x^2-z^2}}{\mathrm{arccosh}\frac{x}{\sqrt{x^2-z^2}}}\Big),\\
\tau=\frac{\log \sqrt{x^2-z^2}}{\sin v}, ~ \text{where}  ~ -\pi/2\le v \le \pi/2, ~ \tau \in \bR^+.
\end{gathered} \tag{3.7}
\end{equation}
\item $y=0$, $z\ne 0$ and $x^2-z^2 = 1$;
\begin{equation}
\begin{gathered}
u=\frac{\pi}{2}, ~v=0, ~ \tau=\mathrm{arccosh}(x), ~ \tau \in \bR^+.
\end{gathered} \tag{3.8}
\end{equation}
\item $y, z =0$;
\begin{equation}
u=0,~v=\frac{\pi}{2},~\tau=\log (x), ~ \tau \in \bR^+. \tag{3.9} 
\end{equation}
\end{enumerate} ~ ~ $\square$
\end{lem}
It is clear that $Y \in \cA\cS^X_{P_1P_2}(\lambda) ~ \text{iff} ~ d^X(P_1,Y)=\lambda \cdot d^X(Y,P_2)$ and $d(Y,P_2)=d(Y^{\mathcal{F}},P_2^{\mathcal{F}})$, where $\mathcal{F}$ 
is a composition of isometries of geometry $X$ which maps $Y$ onto $A_0(1,1,0,0)$ ($ X \in \{ \SXR, \HXR \}$). Therefore, apply the Lemmas 3.1-2 the length 
of a given geodesic segment (i.e. the distance between the two points) and thus it is comparable to $d^X(P_1,Y)$.

This method clearly leads to the following implicit equation of the Apollonius surfaces $\cA\cS^X_{P_1P_2}(\lambda)$ of two proper points 
$P_1(1,a,b,c)$ and $P_2(1,d,e,f)$ with given ratio $\lambda \in \bR^+_0$, in $X$ geometry:
\begin{Theorem}
The implicit equation of the Apollonius surfaces $\cA\cS^X_{P_1P_2}(\lambda)$ of two proper points $P_1(1,a,b,c)$ and $P_2(1,d,e,f)$ with given ratio 
$\lambda \in \bR^+_0$, in $X$ geometry:
\begin{equation}\label{heqeq2}
\begin{gathered}
\cA\cS^X_{P_1P_2}(\lambda)(x,y,z) \Rightarrow \\
4 \omega_X^2{\left(\frac{ax \pm by \pm cz}{\sqrt{a^2 \pm b^2 \pm c^2}\sqrt{x^2 \pm y^2 \pm z^2}} \right)}+{\log^2{\Big(\frac{a^2 \pm b^2 \pm c^2}{x^2 \pm y^2 \pm z^2}\Big)}}=\\
=\lambda^2 \Big[4 \omega_X^2{\left(\frac{dx \pm ey \pm fz}{\sqrt{d^2 \pm e^2 \pm f^2}\sqrt{x^2 \pm y^2 \pm z^2}} \right)}+{\log^2{\Big(\frac{d^2 \pm e^2 \pm f^2}
{x^2 \pm y^2 \pm z^2}\Big)}}\Big],
\end{gathered} \notag
\end{equation}
where if $X=\SXR$ then all $\pm$ signs are $+$, $\omega_X(x)=\arccos(x)$ and if $X=\HXR$ then the all $\pm$ signs are $-$, $\omega_X(x)=\arccosh(x)$. ~~$\square$ 
\end{Theorem}
\begin{figure}[ht]
\centering
\includegraphics[width=13cm]{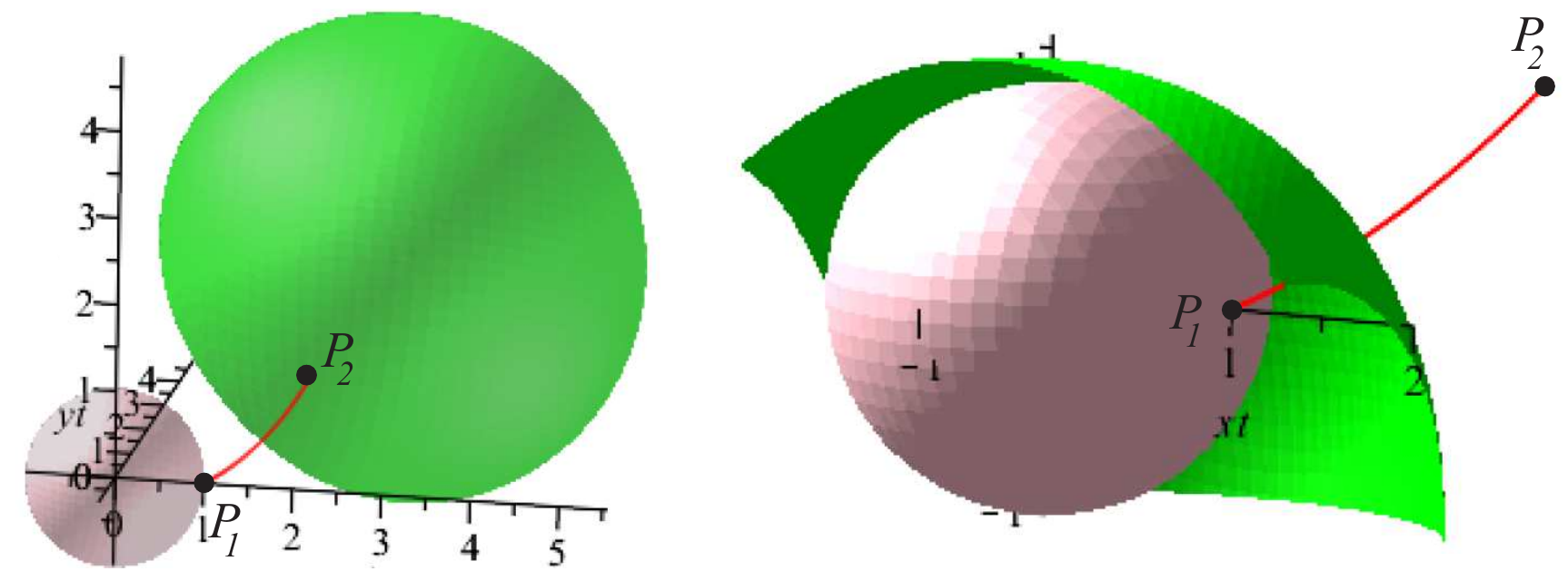}
\caption{The Apollonius surface $\mathcal{A}\cS^{\SXR}_{P_1P_2}(\lambda)$ where $P_1=(1,1,0,0)$, $P_2=(1,2,1,1)$, $\lambda=2$ (left) and $\lambda=1$ (right, equidistance surface)}
\label{}
\end{figure}
\begin{figure}[ht]
\centering
\includegraphics[width=13cm]{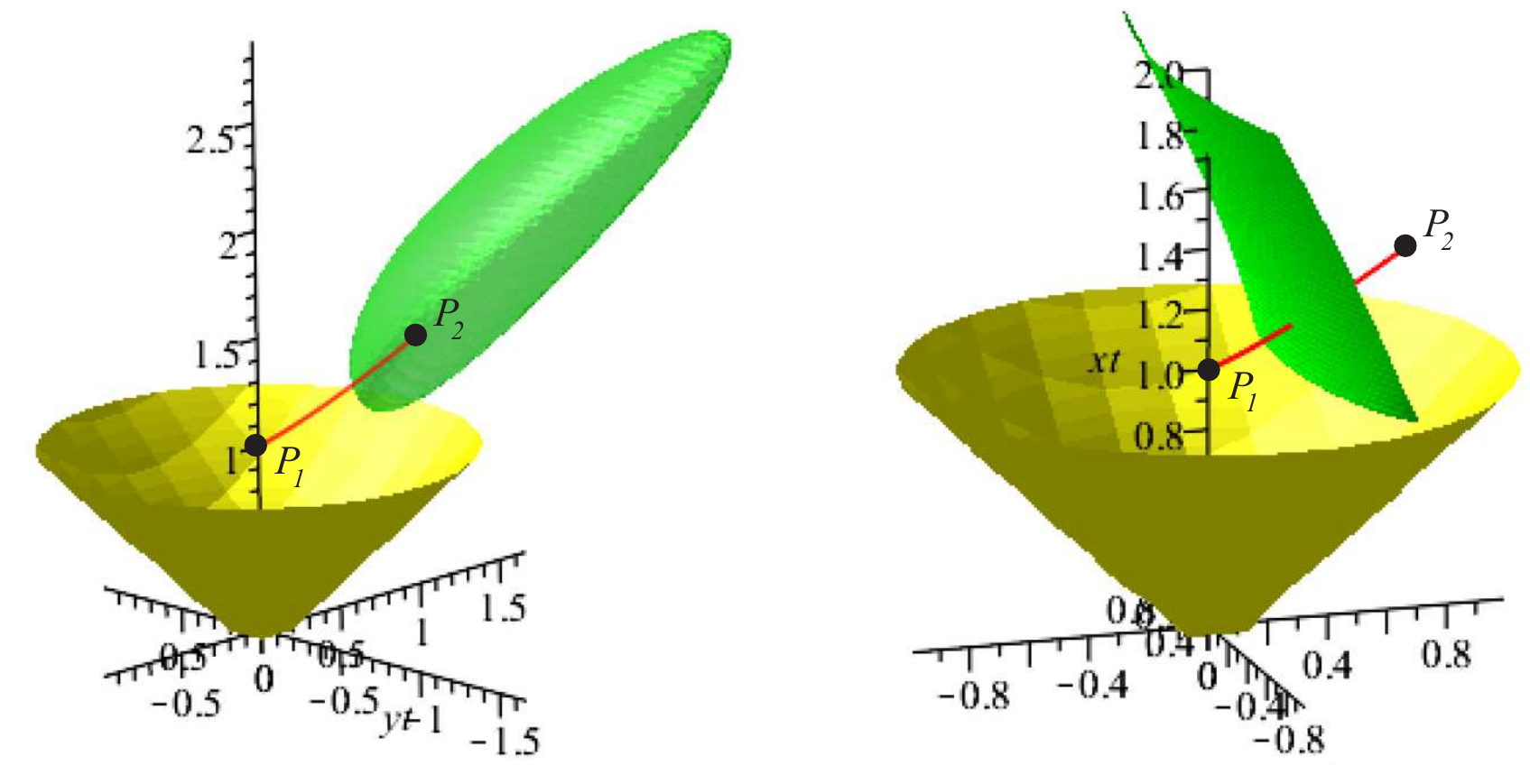}
\caption{The Apollonius surface $\mathcal{A}\cS^{\HXR}_{P_1P_2}(\lambda)$ where $P_1=(1,1,0,0)$, $P_2=(1,3/2,1/1,-1/2)$, $\lambda=2$ (left) and $\lambda=1$ (right, equidistance surface)}
\label{}
\end{figure}

We will use the statements of the following lemma
\begin{lem}[J.~Sz, \cite{Sz20}]
Let $P$ be an arbitrary point and $g^{X}(P_1,P)$ $(X \in \{ \SXR, \HXR \}$, $P_1=(1,1,0,0)$) is a geodesic curve in the considered model of $X$ geometry. 
The points of the geodesic curve $g^X(P_1,P)$ and the centre of the model $E_0=(1,0,0,0)$ lie in a plane in Euclidean sense. ~ ~ $\square$
\end{lem}
\subsection{The surfaces of geodesic triangles}
We consider $3$ points $A_0$, $A_1$, $A_2$ in the projective model of $X$ space (see Section 2) $(X\in\{\SXR, \HXR \})$.
The {\it geodesic segments} $a_k$ connecting the points $A_i$ and $A_j$
$(i<j,~i,j,k \in \{0,1,2\}, k \ne i,j$) are called sides of the {\it geodesic triangle} with vertices $A_0$, $A_1$, $A_2$.

However, defining the surface of a geodetic triangle in $X$ space is not straightforward. The usual geodesic triangle surface definition in these geometries is 
not possible because the geodesic curves starting from different vertices and ending at points of the corresponding opposite edges define different 
surfaces, i.e. {\it geodesics starting from different vertices and ending at points on the corresponding opposite side usually do not intersect 
(Fig.~3 illustrates it in $\SXR$ space)}.
\begin{figure}[ht]
\centering
\includegraphics[width=12cm]{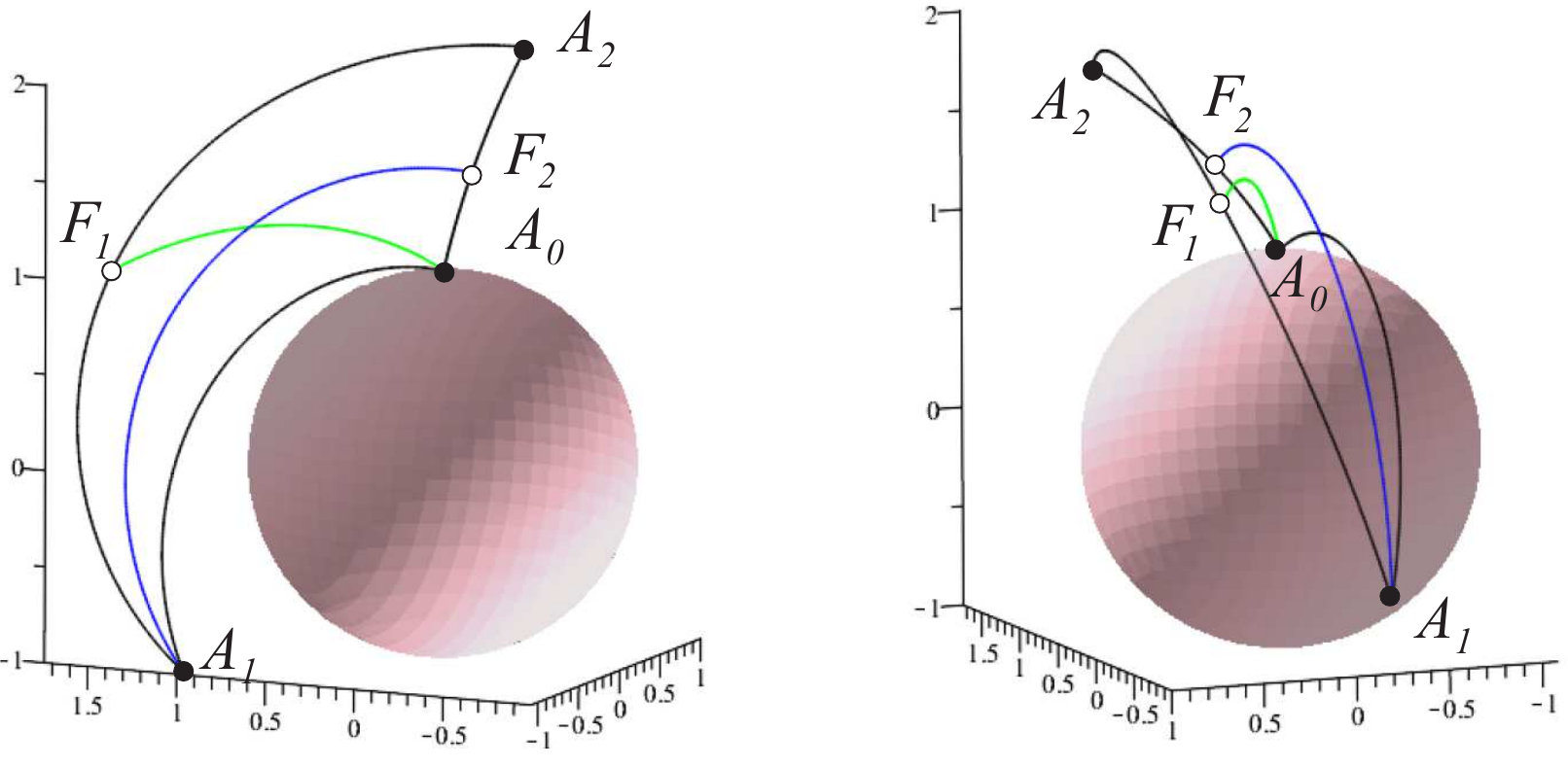}
\caption{$A_0A_1A_2$ is a geodesic triangle in $\SXR$ space with vertices $A_0=(1,1,0,0)$, $A_1=(1,-1,-1,1)$, $A_2=(1,2,1,0)$ and  $F_1$ is the midpoint of 
the geodesic segment $A_1A_2$, $F_2$ is the midpont of geodesic segment $A_0A_2$ The right figure shows that the geodesic segments $A_1F_2$ and $A_0F_1$ do not intersect each other. 
}
\label{}
\end{figure}
Therefore, we introduce a new definition of the surface $\mathcal{S}_{A_0A_1A_2}$ of the geodesic triangle by the following steps:
\begin{definition}
\begin{enumerate}
\item We consider the geodesic triangle $A_0A_1A_2$ in the projective model of $X$ space $(X\in\{\SXR, \HXR \})$ and consider the Apollonius surfaces 
$\mathcal{A}\cS^{X}_{A_0A_1}$ $(\lambda_1)$ and $\mathcal{A}\cS^{X}_{A_2A_0}(\lambda_2)$ ($\lambda_1,\lambda_2 \in [0,\infty)$, $\lambda_1^2+\lambda_2^2>0$). It is clear, that if 
$Y \in \mathcal{C}(\lambda_1,\lambda_2):=\mathcal{A}\cS^{X}_{A_0A_1}(\lambda_1)\cap \mathcal{A}\cS^{X}_{A_2A_0}(\lambda_2)$ then 
$\frac{d^X(A_0,Y)}{d^X(Y,A_1)}=\lambda_1$ and  $\frac{d^X(A_2,Y)}{d^X(Y,A_0)}=\lambda_2$ $\Rightarrow$ $\frac{d^X(A_2,Y)}{d^X(Y,A_1)}=\lambda_1 \cdot \lambda_2$
for parameters $\lambda_1,\lambda_2 \in (0,\infty)$ and if $\lambda_1=0$ then $\mathcal{C}(\lambda_1,\lambda_2)=A_0$, if $\lambda_2=0$ then $\mathcal{C}(\lambda_1,\lambda_2)=A_2$ 
\item 
\begin{equation}
\begin{gathered}
P^X(\lambda_1,\lambda_2):=\{ P \in X~|~P \in \mathcal{C}(\lambda_1,\lambda_2)~{\text{and}}~ d^X(P,A_0)=\min_{Q \in \mathcal{C}(\lambda_1,\lambda_2)}({d^X(Q,A_0)}) \\
~\text{with given real parameters}~\lambda_1,\lambda_2 \in [0,\infty),~ \lambda_1^2 + \lambda_2^2 > 0 \}
\end{gathered} \tag{3.10}
\end{equation}
\item The surface $\mathcal{S}^X_{A_0A_1A_2}$ of the geodesic triangle $A_0A_1A_2$ is 
\begin{equation}
\mathcal{S}_{A_0A_1A_2}^X:=\{P^X(\lambda_1,\lambda_2) \in X,~ \text{where}~\lambda_1,\lambda_2 \in [0,\infty),~ \lambda_1^2 + \lambda_2^2 > 0 \}. \tag{3.11}
\end{equation}
\end{enumerate}
\end{definition}

We introduce the following notations: 

{\it 
1. If the surface of the geodesic triangle is a plane in Euclidean sense which contains the centre of the model $E_0=(1,0,0,0)$, 
then it is called fibre type triangle. 

2. In other cases the triangle is in general type.}
\begin{definition}
Let $\mathcal{S}^X_{A_0A_1A_2}$ be the surface of the geodesic triangle $A_0A_1A_2$ and $P_1$, $P_2 \in \mathcal{S}^X_{A_0A_1A_2}$ two given point. 

1. If the geodesic triangle is in fibre type, then the connecting curve 
$\mathcal{G}^{\mathcal{S}^X_{A_0A_1A_2}}_{P_1P_2} \subset \mathcal{S}^X_{A_0A_1A_2}$ is a geodesic curve in $X$ space. 

2. In other cases the connecting curve $\mathcal{G}^{\mathcal{S}^X_{A_0A_1A_2}}_{P_1P_2}$ is the image of the geodesic curve $g^X_{P_1P_2}$ into the surface $\mathcal{S}^X_{A_0A_1A_2}$
by a central projection with centre $E_0=(1,0,0,0)$.
\end{definition}
\begin{figure}[ht]
\centering
\includegraphics[width=11cm]{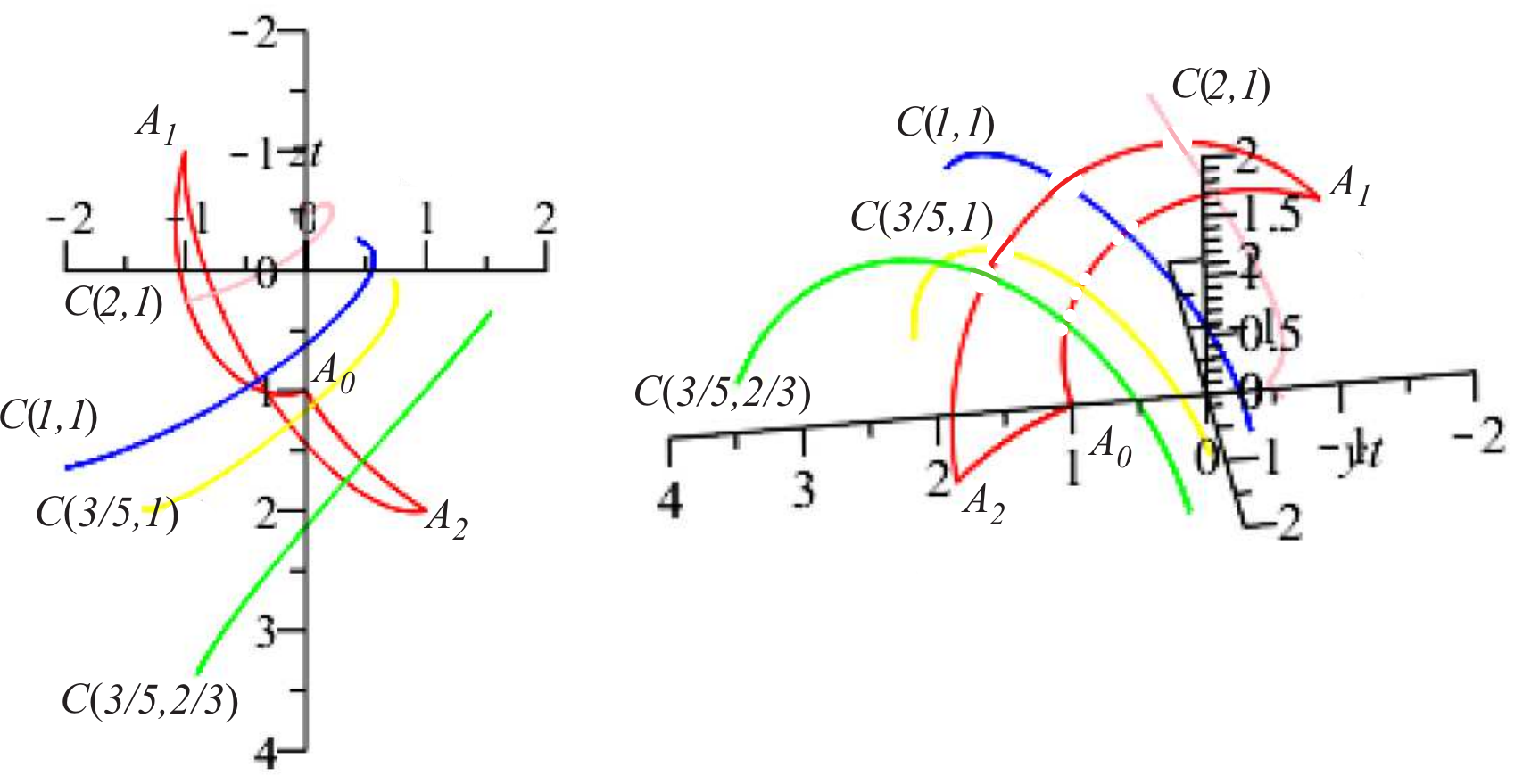}
\caption{Intersection curves $\mathcal{C}(\lambda_1,\lambda_2):=\mathcal{A}\cS^{\SXR}_{A_0A_1}(\lambda_1)\cap \mathcal{A}\cS^{\SXR}_{A_2A_0}(\lambda_2)$ 
related to geodesic triangle $A_0A_1A_2$ in $\SXR$ geometry. 
}
\label{}
\end{figure}
\begin{figure}[ht]
\centering
\includegraphics[width=11cm]{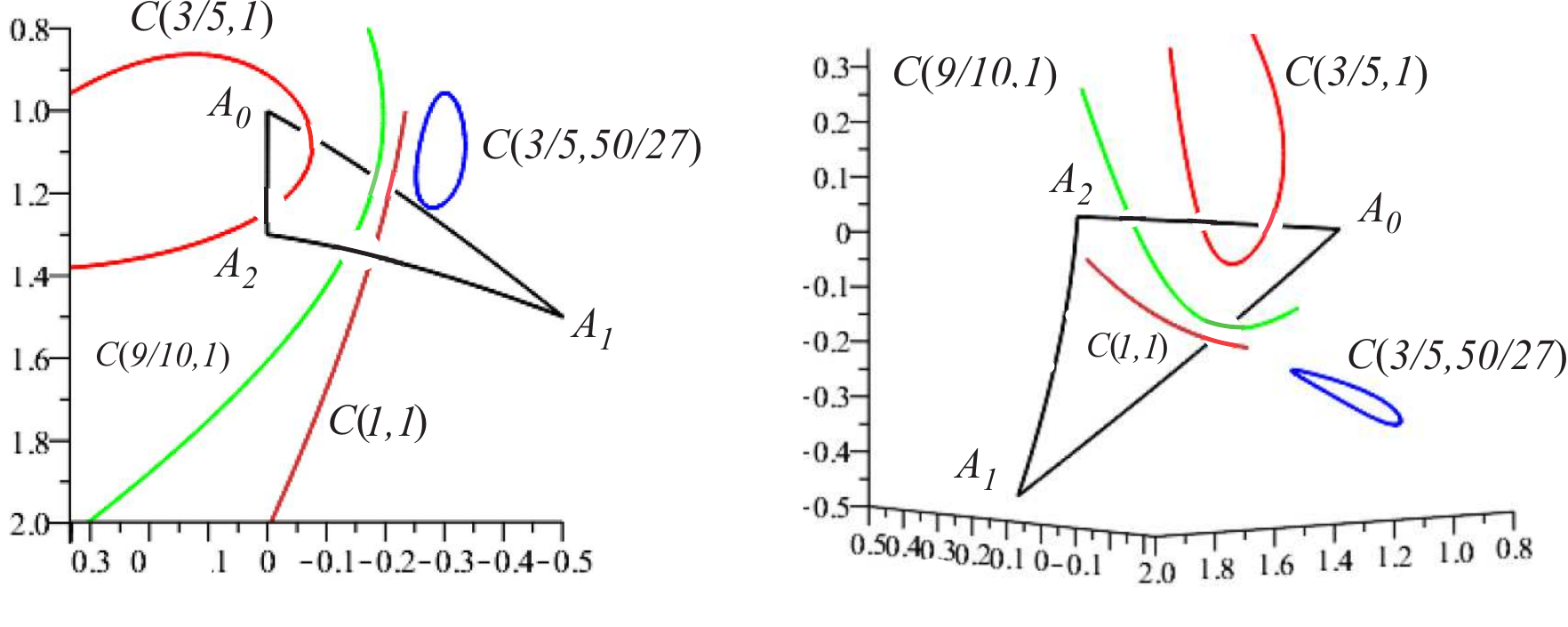}
\caption{Intersection curves $\mathcal{C}(\lambda_1,\lambda_2)=\mathcal{A}\cS^{\HXR}_{A_0A_1}(\lambda_1)\cap \mathcal{A}\cS^{\HXR}_{A_2A_0}(\lambda_2)$ 
related to geodesic triangle $A_0A_1A_2$ in $\HXR$ geometry.
}
\label{}
\end{figure}
\subsection{Geodesic tetrahedra and their circumscribed spheres}
We consider $4$ points $A_0$, $A_1$, $A_2$, $A_3$ in 
the projective model of $X$ space (see Section 2, $X\in\{\SXR, \HXR \})$.
These points are the vertices of a {\it geodesic tetrahedron} in the $X$ space if any two {\it geodesic segments} connecting the points $A_i$ and $A_j$
$(i<j,~i,j \in \{0,1,2,3\})$ do not have common inner points and any three vertices do not lie in a same geodesic curve.
Now, the geodesic segments $A_iA_j$ are called edges of the geodesic tetrahedron $A_0A_1A_2A_3$.
\begin{figure}[ht]
\centering
\includegraphics[width=11cm]{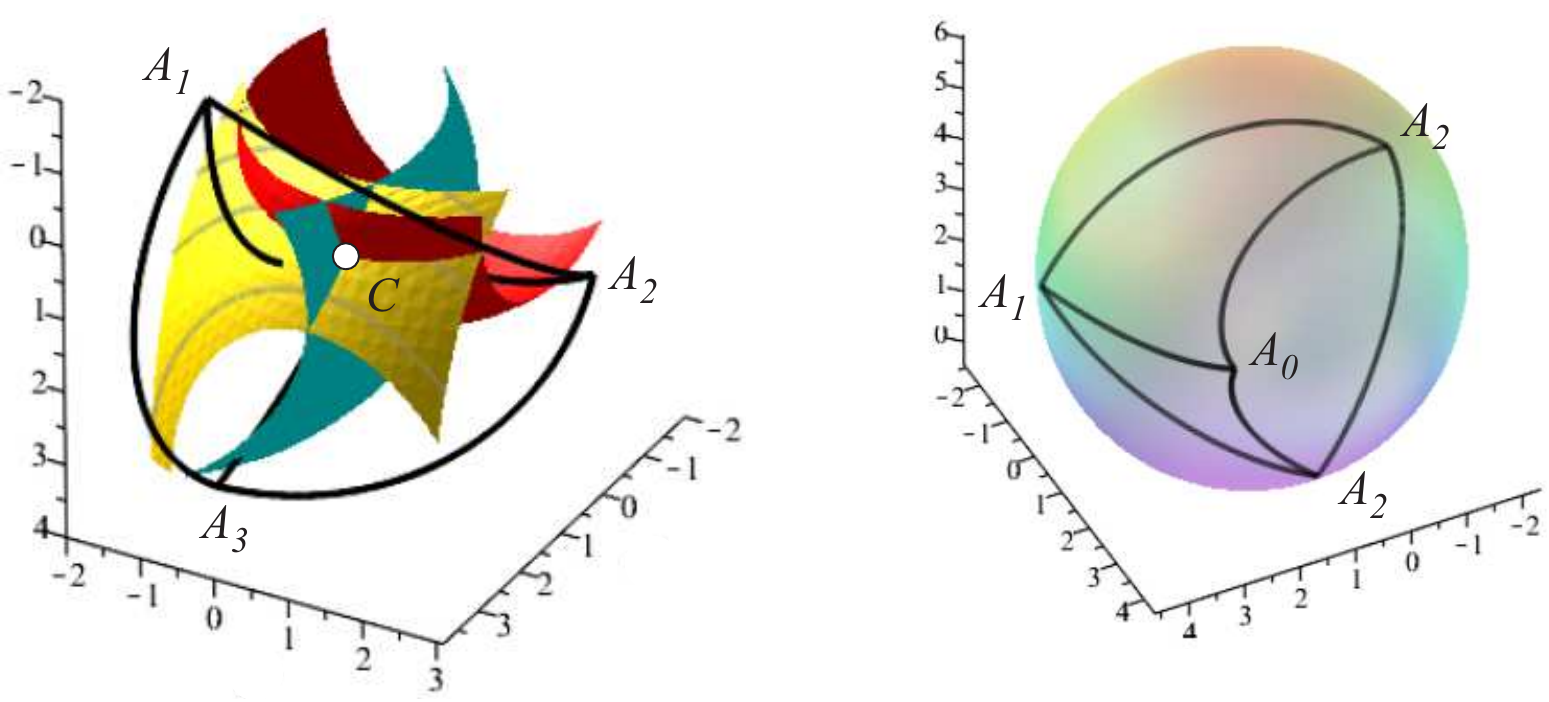}
\caption{Geodesic $\SXR$ tetrahedron with vertices $A_0=(1,1,0,0)$, $A_1=(1,-2, -1/2,3)$, $A_2=(1,1,3,0)$, $A_3=(1,4,-1,2)$
and its circumscibed sphere of radius $r \approx 1.30678$ with circumcenter $C=(1,\approx 0.64697, \approx 0.51402, \approx 0.15171)$.}
\label{}
\end{figure}
The circumscribed sphere of a geodesic tetrahedron is a geodesic sphere (see Definitions 2.2, 2.6, and Fig.~6,7,8,9) that touches each of the tetrahedron's vertices.
As in the Euclidean case the radius
of a geodesic sphere circumscribed around a tetrahedron $T$ is called the circumradius of $T$, and the center point of this sphere is called the circumcenter of $T$.
\begin{figure}[ht]
\centering
\includegraphics[width=10cm]{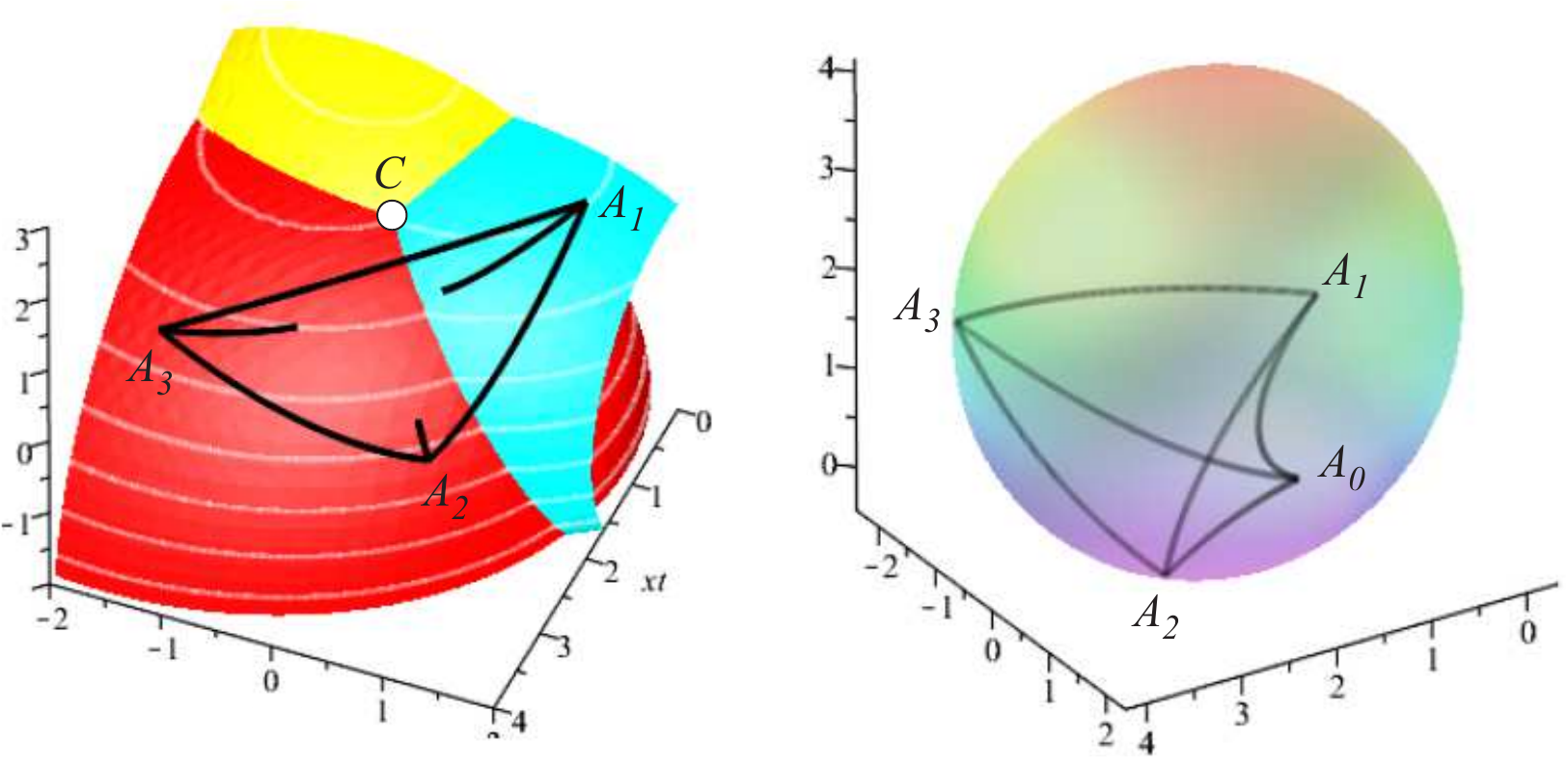}
\caption{Geodesic $\SXR$ tetrahedron with vertices $A_0=(1,1,0,0)$, $A_1=(1,2,2,3)$, $A_2=(1,3,1,0)$, $A_3=(1,4,-1,2)$
and its circumscibed sphere of radius $r \approx 0.97720$ with circumcenter $C=(1,\approx 1.34000, \approx -0.01705, \approx 1.27323)$. }
\label{}
\end{figure}

\subsubsection{Circumscribed spheres in $\SXR$ space}

The next Lemma follows directly from the properties of the geodesic distance function of $\SXR$ space (see Definition 2.1 and (2.14)):
\begin{lem}
For any $\SXR$ geodesic tetrahedron there exists uniquely a geodesic surface on which all four vertices lie. If its radius less or equal to $\pi$ then 
the above surface is a geodesic sphere (called circumscibed sphere, see Definition 2.2). \ \ $\square$
\end{lem}
The procedure to determine the radius and the circumcenter of a given geodesic $\SXR$ tetrahedron is the folowing:

The circumcenter $C=(1,x,y,z)$ of a given translation tetrahedron $A_0A_1A_2A_3$ $(A_i=(1,x^i,y^i,z^i), ~ i \in \{0,1,2,3\})$
has to hold the following system of equation:
\begin{equation}
d^{\SXR}(A_0,C)=d^{\SXR}(A_1,C)=d^{\SXR}(A_2,C)=d^{\SXR}(A_3,C), \tag{3.12}
\end{equation}
therefore it lies on the geodesic-like bisector surfaces $\cA\cS^{\SXR}_{A_i,A_j}$ $(i<j,~i,j \in \{0,1,2,3\}$) which equations are determined in Lemma 3.3.
The coordinates $x,y,z$ of the circumcenter of the circumscribed sphere around the tetrahedron $A_0A_1A_2A_3$ are obtained by the system of equation
derived from the facts:
\begin{equation}
C \in \cA\cS^{\SXR}_{A_0A_1}(1), \cA\cS^{\SXR}_{A_0A_2}(1), \cA\cS^{\SXR}_{A_0A_3}(1). \tag{3.13}
\end{equation}
Finally, we get the circumradius $r$ as the geodesic distance e.g. $r=d^{\SXR}(A_1,C)$.

We applied the above procedure to two tetrahedra determined their centres and the
radii of their circumscribed balls that are described in  Fig.~6 and 7.~ $\square$
\begin{figure}[ht]
\centering
\includegraphics[width=10cm]{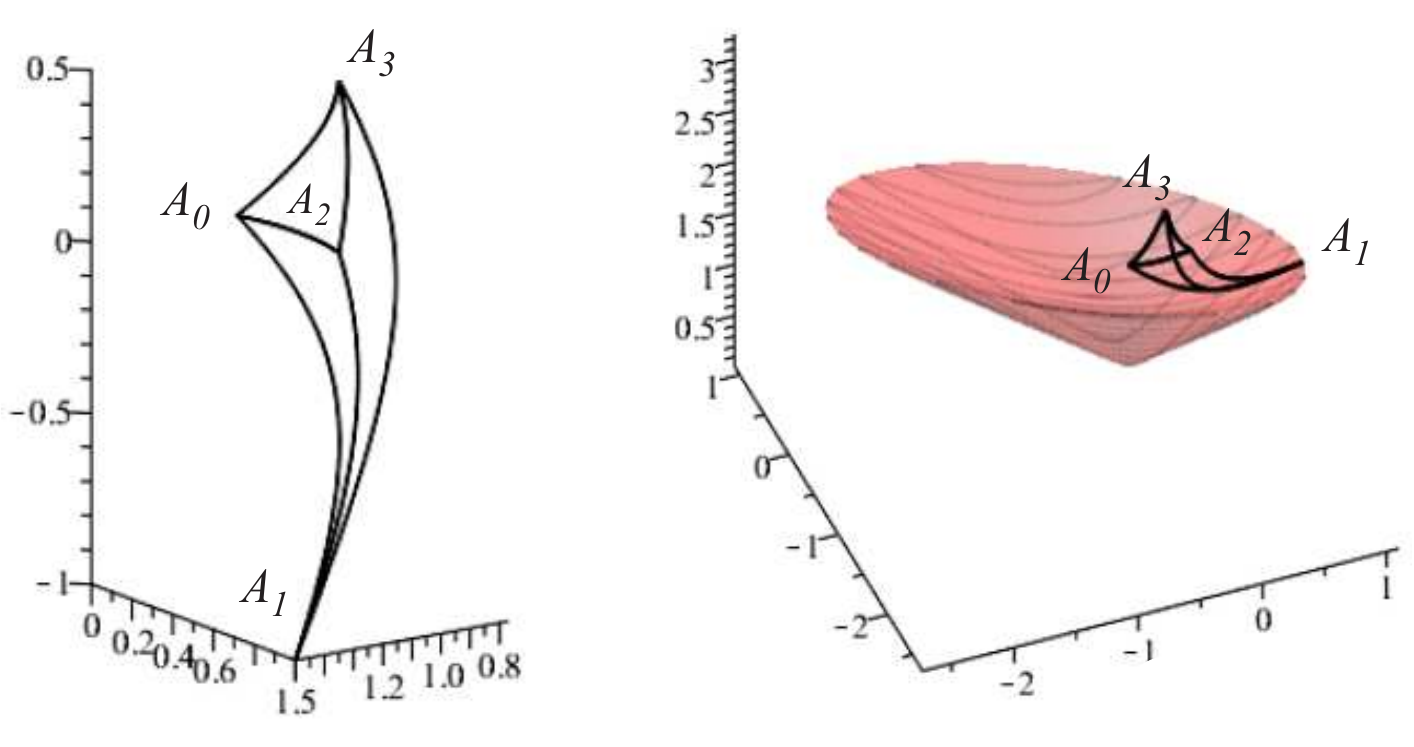}
\caption{Geodesic $\HXR$ tetrahedron with vertices $A_0=(1,1,0,0)$, $A_1=(1,3/2,1,-1)$, $A_2=(1,1,1/2,0)$, $A_3=(1,1,1/2,1/2)$
and its circumscibed sphere of radius $r \approx 2.89269$ with circumcenter $C=(1,\approx 0.07017, \approx -0.02714, \approx -0.02640)$. }
\label{}
\end{figure}
\begin{figure}[ht]
\centering
\includegraphics[width=10cm]{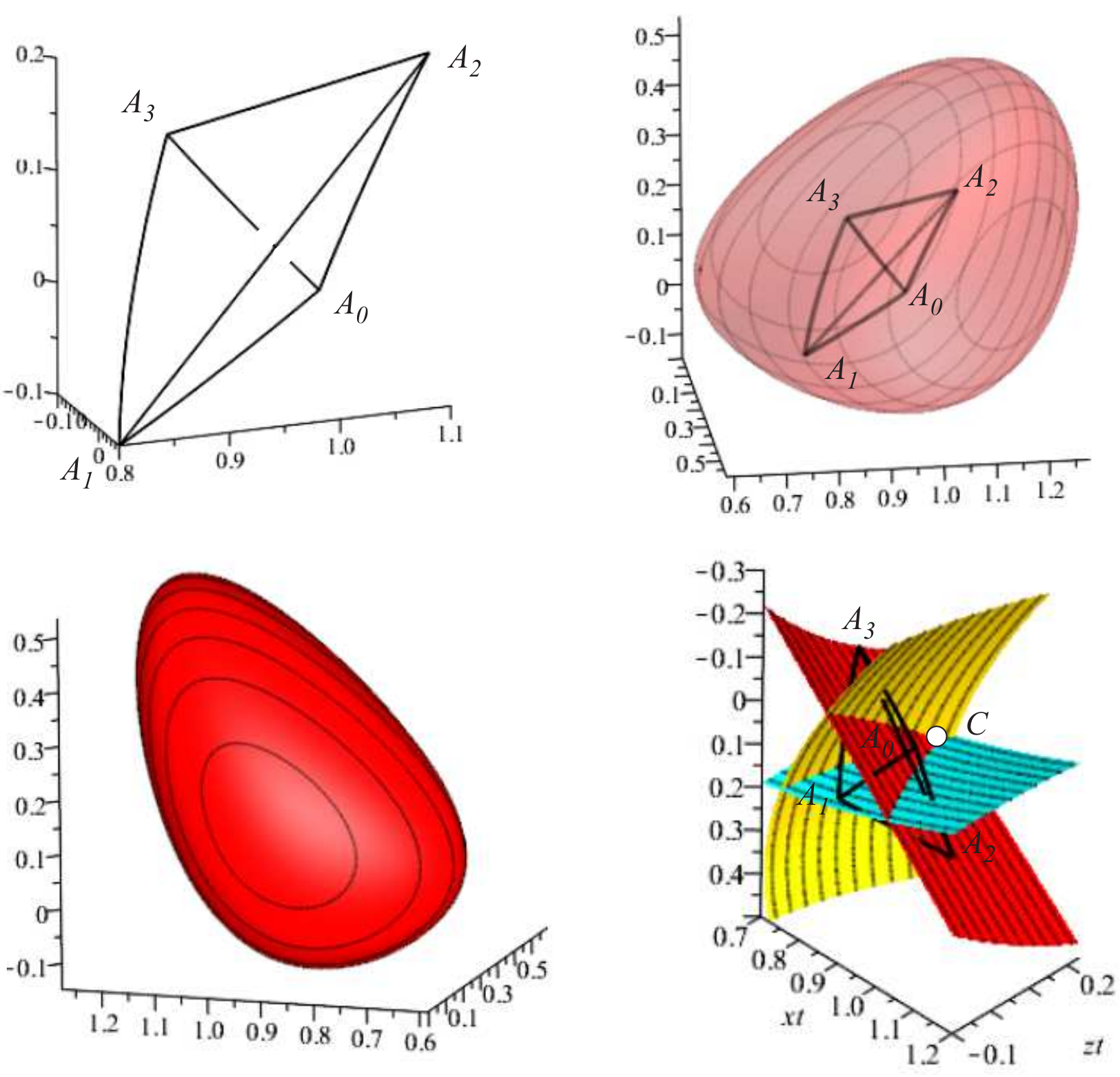}
\caption{Geodesic $\HXR$ tetrahedron with vertices $A_0=(1,1,0,0)$, $A_1=(1,9/10,3/25,-1/10)$, $A_2=(1,11/10,1/5,0)$, $A_3=(1,4/5,-1/10,1/20)$
and its circumscibed sphere of radius $r \approx 0.37162$ with circumcenter $C=(1,\approx 0.85902, \approx 0.16307, \approx 0.20314)$. }
\label{}
\end{figure}
\subsubsection{Circumscribed spheres in $\HXR$ space}
In the $\HXR$ geometry the procedure to determine the radius and the circumcenter of a given geodesic $\HXR$ tetrahedron is similar to the $\SXR$ case but not the same.
The circumcenter $C$ and its circumradius $r$ of the circumscribed sphere around the tetrahedron $A_0A_1A_2A_3$ is obtained by the system of equation
derived from the bisector surfaces of the edges of the given tetrahedron:
\begin{equation}
C \in \cA\cS^{\HXR}_{A_0A_1}(1), \cA\cS^{\HXR}_{A_0A_2}(1), \cA\cS^{\HXR}_{A_0A_3}(1), \tag{3.14}
\end{equation}
but the above surfaces do not always intersect at a proper point of $\HXR$ space. It can be an outer point (relating to the model) or a point at the infinity (a point on the boundary of the model), 
similarly to the hyperbolic spaces.
\begin{rmrk}
If the common point of the above bisectors lies at the infinity then the vertices of tetrahedron lie on a horosphere-like surface and if the 
common point is an outer point then the vertices of the tetrahedron are on a hypershere-like surface. 
These surfaces will be examined in detail in a forthcoming paper.
\end{rmrk}
The next Lemma follows directly from the properties of the geodesic distance function of $\HXR$ space (see Definition 2.6 and (2.10)):
\begin{lem}
For any $\HXR$ geodesic tetrahedron there exists uniquely a geodesic surface on which all four vertices lie. If its centre is a proper point of $\HXR$ space then 
the above surface is a geodesic sphere (called circumscibed sphere, see Definition 2.6). \ \ $\square$
\end{lem}

We applied the above procedure to two tetrahedra determined their centres and the
radii of their circumscribed balls that are described in  Fig.~8 and 9.~ $\square$
\section{Menelaus' and Ceva's theorems in $\SXR$ and $\HXR$ spaces} 
First we recall the definition of simply ratios in spherical $\bS^2$ and the hyperbolic $\bH^2$ planes. 
The models of the above plane geometries of constant curvature are embedded in the models of the previously described geometries $\SXR$ and $\HXR$ 
as ``base planes" and are used hereinafter for our discussions. 

A spherical triangle is the space included by arcs of great circles on the surface of a sphere, subject to the limitation that these arcs and the 
further circle arcs in the following items in the spherical plane are always less or equal than a semicircle.
\begin{definition}
If $A$, $B$ and $P$ are distinct points on a line in the $Y\in\{\bH^2,\bS^2\}$ space, then 
their simply ratio is
$s^Y(A,P,B) =  w^Y(d^Y(A,P))/w^Y(d^Y(P,B))$, if $P$ is between $A$ and $B$, and
$s^Y(A,P,B) = -w^Y(d^Y(A,P))/w^Y(d^Y(P,B))$, otherwise where
$w^Y(x):=sin(x)$ if $Y=\bS^2$ and $w^Y(x):=sinh(x)$ if $Y=\bH^2$. 
\end{definition}

\begin{rmrk} Basic properties of simply ratio:
\begin{enumerate}
\item $s^Y(A,P,B) = -s^Y(B,P,A)$,
\item if $P$ is between $A$ and $B$, then $s^Y(A,P,B) \in (0,1)$,
\item if $P$ is on $AB$, beyond $B$, then $s^Y(A,P,B) \in (-\infty,-1)$,
\item if $P$ is on $AB$, beyond $A$, then $s^Y(A,P,B) \in (-1,0)$.
\end{enumerate}
\end{rmrk}
Note that the value of $s^Y(A,P,B)$ determines the position of $Y$ relative to $A$ and $B$.

With this definition, the corresponding sine rule of geometry $Y$ leads to Menelaus's and Ceva's theorems \cite{K,PS14}:

\begin{Theorem} [Menelaus's Theorem for triangles in $Y$ plane]
If is a $l$ line not through any vertex of an triangle $ABC$ such that
$l$ meets $BC$ in $Q$, $AC$ in $R$, and $AB$ in $P$,
then $$s^Y(A,P,B)s^Y(B,Q,C)s^Y(C,R,A) = -1.$$ ~ ~ $\square$
\end{Theorem}

\begin{Theorem}[Ceva's Theorem for triangles in $Y$ plane]
If $T$ is a point not on any side of a triangle $ABC$ such that
$AT$ and $BC$ meet in $Q$, $BX$ and $AC$ in $R$, and $CX$ and $AB$ in $P$,
then $$s^Y(A,P,B)s^Y(B,Q,C)s^Y(C,R,A) = 1.$$ ~ ~ $\square$
\end{Theorem}
\subsection{Generalizations of Menelaus' and Ceva's theorems}
\subsubsection{Geodesic triangle in general position}
First we consider a {\it general location geodesic triangle} $A_0A_1A_2$ in the projective model of $X$ space (see Section 3.1) $(X\in\{\SXR, \HXR \})$.
Without limiting generality, we can assume that $A_0=(1,1,0,0)$ and $A_2$ lies in the coordinate plane $[x,y]$. 
The geodesic lines that contain the sides $A_0A_1$ and $A_0A_2$ of the given triangle can be characterized directly by the corresponding parameters $v$ and $u$ (see (2.4) and (2.8)).

The geodesic curve including the side segment $A_1A_2$ is also determined by one of its endpoints and its parameters but in order to determine the corresponding parameters of this 
geodesic line we use {\it orientation preserving isometric transformations} $\bT^{X}(A_2)$, as elements of the isometry group of $X$ geometry, that
maps the $A_2=(1,x_2,y_2,0)$ onto $A_0=(1,1,0,0)$ (up to a positive determinant factor). 

\begin{rmrk}
I note here, that this orientation preserving isometry is not unique, but in all cases the $v$ parameters of the 
``image geodesics" are equal (of course the $u$ parameters may be different) and in the further derivation, only the values of the parameter $v$ will be needed.
\end{rmrk} 
We extend the definition of the simply ratio to the $X\in\{\SXR, \HXR \}$ spaces.
If $X=\SXR$ then is clear that the space contains its ``base sphere" (unit sphere centred in $E_0$) which is a geodesic surface. 
Therefore, similarly to the spherical spaces we assume that the geodesic arcs  
in the following items are always less or equal than a semicircle.
\begin{definition}
If $A$, $B$ and $P$ be distinct points on a non-fibrum-like geodetic curve in the $X \in \{\SXR, \HXR \}$ space , then 
their simply ratio is
$$s_g^X(A,P,B) =  w^X\Big({d^X(A,P)}{\cos(v)}\Big)/w^X\Big({d^X(P,B)}{\cos(v)}\Big),$$ if $P$ is between $A$ and $B$, and
$$s_g^X(A,P,B) = -w^X\Big({d^X(A,P)}{\cos(v)}\Big)/w^X\Big({d^X(P,B)}{\cos(v)}\Big),$$ 
otherwise where
$w^X(x):=sin(x)$ if $X=\SXR$, $w^Y(x):=sinh(x)$ if $X=\HXR$ and $v$ is the parameter of the geodesic curve containing points $A,B,P$ (see Fig.~10).  
\end{definition}
\begin{Theorem}[Ceva's Theorem for triangles in general location]
If $T$ is a point not on any side of a geodesic triangle $A_0A_1A_2$ in $X\in\{\SXR, \HXR\}$ such that
the curves $A_0T$ and $g_{A_1A_2}^X$ meet in $Q$, $A_1T$ and $g_{A_0A_2}^X$ in $R$, and $A_2T$ and $g_{A_0A_1}^X$ in $P$, $(A_0T, A_1T, A_2T \subset \mathcal{S}^X_{A_0A_1A_2})$
then $$s_g^X(A_0,P,A_1)s_g^X(A_1,Q,A_2)s_g^X(A_2,R,A_0) = 1.$$ 
\end{Theorem}
{\bf{Proof:}}

\begin{figure}[ht]
\centering
\includegraphics[width=13cm]{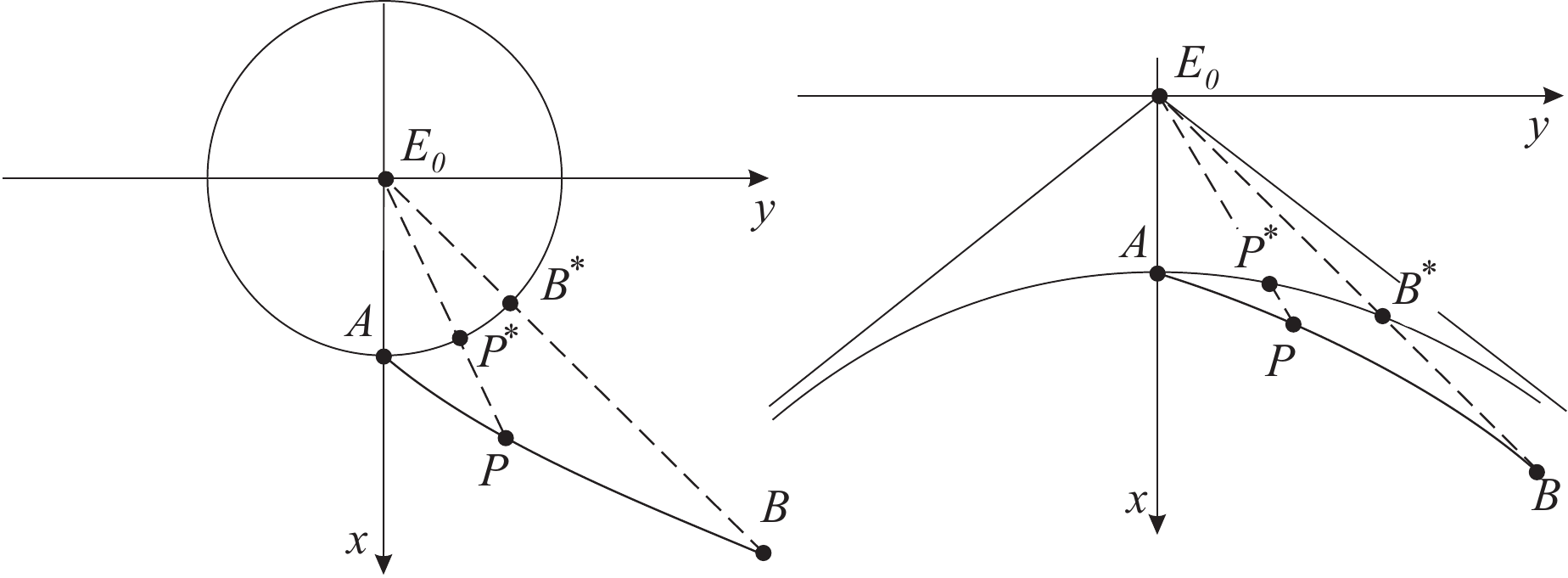}
\caption{Simply ratios in $\SXR$ and $\HXR$ spaces}
\label{}
\end{figure}
Let $A$, $B$ and $P$ be distinct points on a non-fibrum-like geodetic curve in the $X$ and let $A^*$, $B^*$ and $P^*$ their centrally projected images from 
the origin $E_0$ into the base plane of $X$ (see Fig.~10). We assume without loss of genarality, that $A=A^*=(1,1,0,0)$ and $B$ lies in the $[x,y]$ coordinate plan.
Therefore, the points of the geodesic curve segment $g_{AB}^X$ also are included in it (see Lemma 3.3). Let the parameters of $g_{AB}^X$ be $(t,u=0,v)$ (see (2.4), (2.8)) 
and the geodesic curve $g_{AP}^X$
are $(t_p,u_p=0,v_p=v)$. Their images $g_{A^*P^*}^X$ and $g_{A^*B^*}^X$ are determined by parameters $(t_p\cdot \cos{v}, u_p=0, v_p=0)$ and $(t\cdot \cos{v}, u=0, v=0)$:
\begin{equation}
d^X(A^*,P^*)=d^X(A,P)\cdot \cos{v},\ \ \ d^X(P^*,B^*)=d^X(P,B)\cdot \cos{v}. \notag
\end{equation}
Let $A_0^*A_1^*A_2^*$ $(A_0^*=A_0)$ be the centrally projected image of the geodesic triangle $A_0A_1A_2$ from the origin $E_0$ into the base plane of $X$ (see Fig.~10).

Moreover, let $T$ is a point not on any side of a geodesic triangle $A_0A_1A_2$ such that
the curves $A_0T$ and $g_{A_1A_2}^X$ meet in $Q$, $A_1T$ and $g_{A_0A_2}^X$ in $R$, and $A_2T$ and $g_{A_0A_1}^X$ in $P$, $(A_0T, A_1T, A_2T \subset \mathcal{S}^X_{A_0A_1A_2})$.

The centrally projected images $A_0^*Q^*$, $A_1^*R^*$ and $A_2^*P^*$ of the curves $A_0Q$, $A_1R$ and $A_2P$ from $E_0$ are geodesic curves in the corresponding ``base plane" of $X$.
$T^*$ and $T$ lie in the planes $A_0QE_0$, $A_1RE_0$ and $A_2PE_0$ (see Lemma 3.3). Therefore, this theorem is direct consequence of the Definition 4.6 of the simply ratio and the corresponding 
Ceva's theorem in spherical and hyperbolic planes (Theorem 4.4). ~ ~ $\square$ 
\begin{Theorem}[Menelaus's theorem for triangles in general location]
If is a $l$ line not through any vertex of an geodesic triangle $A_0A_1A_2$ lying in its surface $\mathcal{S}^X_{A_0A_1A_2}$ in $X\in\{\SXR, \HXR\}$ geometry such that
$l$ meets geodesic curves $g_{A_1A_2}^X$ in $Q$, $g_{A_0A_2}^X$ in $R$, and $g_{A_0A_1}^X$ in $P$,
then $$s^X_g(A_0,P,A_1)s^X_g(A_1,Q,A_22)s^X_g(A_2,R,A_0) = -1.$$
\end{Theorem} 
{\bf{Proof:}} Similarly to the above proof, this theorem is follows by the Definition 4.6 of the simply ratio and the corresponding 
Ceva's theorem in spherical and hyperbolic planes (Theorem 4.4). ~ ~ $\square$ 
\subsubsection{Fibre type triangle}
We consider a {\it fibre type geodesic triangle} $A_0A_1A_2$ in the projective model of $X$ space.
Without limiting generality, we can assume that $A_0=(1,1,0,0)$, $A_1=(1,x_1,y_1,0)$ and $A_2=(1,x_2,y_2,0)$ lies in the coordinate plane $[x,y]$. 
The geodesic lines that contain the sides $A_0A_1$ and $A_0A_2$ of the given triangle can be characterized directly by the corresponding parameters 
$v$ and $u=0$ similarly to the above case.

We extend the definition of the simply ratio to the $X\in\{\SXR, \HXR \}$ spaces.

If $X=\SXR$ is clear that the $\SXR$ space contains its ``base sphere" (unit sphere centred in $E_0$) which is a geodesic surface. 
Therefore, similarly to the spherical spaces we assume that the geodesic arcs  
in the following items are always less or equal than a semicircle.
\begin{definition}
If $A$, $B$ and $P$ be distinct points on fibrum-like geodetic curve in the $X \in \{\SXR, \HXR \}$ space , then 
their simply ratio is
$$s^X_f(A,P,B) =  d^X(A,P)/d^X(P,B)$$ if $P$ is between $A$ and $B$, and
$$s^X_f(A,P,B) = -d^X(A,P)/d^X(P,B)$$ (see Fig.~10).  
\end{definition}
\begin{Theorem}[Ceva's Theorem in $X$ geometry for triangles in fibre types]
If $T$ is a point not on any side of a geodesic triangle $A_0A_1A_2$ in $X\in\{\SXR, \HXR\}$ such that
the geodesic curves $g_{A_0T}^X$ and $g_{A_1A_2}^X$ meet in $Q$, $g_{A_1T}^X$ and $g_{A_0A_2}^X$ in $R$, and $g_{A_2T}^X$ and $g_{A_0A_1}^X$ in $P$, 
$(g_{A_0T}^X, g_{A_1T}^X, g^X_{A_2T} \subset \mathcal{S}^X_{A_0A_1A_2}$
then $$s^X_f(A_0,P,A_1)s^X_f(A_1,Q,A_2)s^X_f(A_2,R,A_0) = 1.$$
\end{Theorem}
\begin{Theorem}[Menelaus's theorem in $X$ geometry for triangles in fibre types]
If is a $l$ line not through any vertex of an geodesic triangle $A_0A_1A_2$ lying in its surface $\mathcal{S}^X_{A_0A_1A_2}$ in $X\in\{\SXR, \HXR\}$ geometry such that
$l$ meets geodesic curves $g_{A_1A_2}^X$ in $Q$, $g_{A_0A_2}^X$ in $R$, and $g_{A_0A_1}^X$ in $P$,
then $$s^X_f(A_0,P,A_1)s^X_f(A_1,Q,A_22)s^X_f(A_2,R,A_0) = -1$$.
\end{Theorem} 
{\bf{Proof of both theorems:}}

Let $P_1=(1,a,b,c)$ and $P_2=(1,d,e,f)$ be proper points of $X$ space. Their distance $d^X(P_1,P_2))$ is (see Theorem 3.3):
\begin{equation}\label{heqeq2}
\begin{gathered}
d^X(P_1,P_2))=\sqrt{ \omega_X^2{\left(\frac{ax \pm by \pm cz}{\sqrt{a^2 \pm b^2 \pm c^2}\sqrt{d^2 \pm e^2 \pm f^2}} \right)}+
{\log^2{\Big(\frac{\sqrt{a^2 \pm b^2 \pm c^2}}{\sqrt{d^2 \pm e^2 \pm f^2}}\Big)}}},
\end{gathered} \notag
\end{equation}
where if $X=\SXR$ then all $\pm$ signs are $+$, $\omega_X(x)=\arccos(x)$ and if $X=\HXR$ then the all $\pm$ signs are $-$, $\omega_X(x)=\arccosh(x)$.
Without loss of generality we can assume, that $P_1=A=(1,1,0,0)$ and $P_2=B=(1,d,e,f)$ (see Fig.~10) and $A^*$ and $B^*$ are their 
centrally projected images from the origin $E_0$ into the base plane of $X$. From the structure of geometries $X$ follows, 
that $d^X(P_1,P_2))=\sqrt{(d^X(A^*,B^*))^2+(d^X(B^*,B))^2}$ where 
$(d^X(A^*,B^*))^2 = \omega_X^2{\left(\frac{d}{\sqrt{\sqrt{d^2 \pm e^2 \pm f^2}}} \right)}$ and 
$(d^X(B^*,B))^2= \log^2{\Big(\frac{{1}}{\sqrt{d^2 \pm e^2 \pm f^2}}\Big)}$.

Therefore, each right angled fibre-like triangle in $X$ space can be uniquely assigned to an Euclidean righ angled triangle with the same side lenghts and 
it is clear, that the same is true for any fibre-like geodesic triangle.

Thus, we can formulate similar theorems for the fibre-like geodesic triangle as for the corresponding Euclidean triangles therefore 
the Ceva's and Menelaus' theorems in $X$ geometry follows from the well-known corresponding Euclidean ones. \ \ $\square$

Similar problems in other homogeneous Thurston geometries
represent another huge class of open mathematical problems. For
$\NIL$, $\SOL$, $\SLR$ geometries only very few results are known
\cite{MSz12}, \cite{MSzV}, \cite{CsSz16},  \cite{Sz19}, \cite{MSz}, \cite{Sz12-1}, \cite{Sz13-1}, \cite{Sz13-2},. Detailed studies are the objective of
ongoing research.
\medbreak

\end{document}